\documentclass[12pt]{article}
 \usepackage[german]{babel}
\usepackage{amssymb}
\usepackage{amsmath}
\usepackage{graphicx}
\usepackage{bbm}
\usepackage[right]{eurosym}
\usepackage[left=2cm,right=2cm,top=2cm,bottom=3cm]{geometry}
\def\vs{\vspace}
\def\noi{\noindent}

\def\IN{\mathbb N}
\def\IZ{\mathbb Z}
\def\IR{\mathbb R}

\def\an{\mathrm{an}}
\def\exp{\mathrm{exp}}

\def\ma{\mathcal}

\begin{document}
\begin{center}
{\bf \Large Asymptotics of parameterized exponential integrals given by Brownian motion on globally subanalytic sets}
\end{center}

\centerline{Tobias Kaiser and Julia Ruppert}

\vspace{0.7cm}\noi \footnotesize {{\bf Abstract.} We consider parameterized exponential integrals coming from the time evolution of the probability distribution of Brownian motion on globally subanalytic sets. We establish definability results and asymptotic expansions.

\normalsize
\section*{Introduction}

Integration in the o-minimal setting is a difficult task. One can this observe already from the fact, that the reciprocal function $x\mapsto 1/x$ is semialgebraic whereas its antiderivative, the logarithmic function $x\mapsto \log x$, is not.  
So far there is no general result as formulated as a question in Van den Dries [8, p. 147] or as a conjecture in [14, p. 1904].

But there are several very deep partial results.
Speissegger [21] has shown that the antiderivative of a function definable in an o-minimal structure on the real field is definable in an o-minimal expansion, to be more precise in its Pfaffian closure.
Comte, Lion and Rolin [5, 20] have studied parameterized integrals $x\mapsto \int f(x,y)\;dy$ where $f(x,y)$ is globally subanalytic. They have shown that the set of parameters where the integral exists is again globally subanalytic and that its evaluation on this set is a finite sum of finite products of globally subanalytic functions and logarithms of positive globally subanalytic functions (which have been later called constructible functions).
In [15], the special case $x\mapsto \int f(x,y)\;dy$ for semialgebraic $f(x,y)$ has been studied.
Cluckers and D. Miller [2, 3, 4] have established that the class of constructible functions is closed unter taking parameterized integrals. In particular, they handle the case $x\mapsto \int f(x,y)\log(g(x,y))\;dy$ for globally subanalytic $f(x,y)$ and $g(x,y)$. More recently, Cluckers et al. [1] have considered oscillatory integrals $x\mapsto\int f(x,y)e^{ig(x,y)}$ for globally subanalytic $f(x,y)$ and $g(x,y)$. These are motivated by Fourier transformation and lead necessarily out of the o-minimal context.
Inside the o-minimal category, the next natural step would be to analyze parameterized integrals of the form $x\mapsto \int f(x,y)e^{g(x,y)}\;dy$ where $f(x,y)$ and $g(x,y)$ are globally subanalytic, the motivation coming from Laplace transformation. But here the situation is completely new. For example, the error function which is up to a multiplicative constant the antiderivative of $x\mapsto e^{-x^2}$ is by Van den Dries et al. [9] not definable in the o-minimal structure $\IR_{\an,\exp}$.

We consider in the present paper parameterized integrals of the form $t\mapsto\int f(y)e^{-y^2/2t}\;dy$ for globally subanalytic or constructible $f(y)$ and single variables $t,y$. We are interested in them because they have a concrete analytic-geometric meaning; they are connected to Brownian motion.
Brownian motion is one of the most important stochastic processes with vast applications (see for example [17]).  Modelling jittery motion it seems at first curious to consider it in the globally subanalytic or more general o-minimal context. On a microscopic level, it is definitely not tame. But on a macroscopical level, one can hope to obtain tameness results when starting in a tame setting.

 \vs{0.5cm}
\hrule

\vs{0.4cm}
{\footnotesize{\itshape 2010 Mathematics Subject Classification:} 03C64, 26A42, 26A12, 32B20, 41A60, 60G15}
\newline
{\footnotesize{\itshape Keywords and phrases:} globally subanalytic sets, o-minimality, exponential parameterized integrals, Brownian motion, asymptotic expanions}

We are not interested in the dynamics of the Brownian motion but in the time evolution of the probability that the Brownian motion $(B_t)_{t\geq 0}$ with start value $a\in \IR^n$ is in a given subset $A$ of $\IR^n$ which is assumed to be globally subanalytic or more 
general definable in an o-minimal structure on the field of reals.
This probabilty at time $t$ is given by
$$P_a(B_t\in A)=
\left\{\begin{array}{ccc}\delta_a(A),&&t=0,\\
&\mbox{if}&\\
\frac{1}{(2\pi t)^{\frac{n}{2}}}\int_A e^{-\frac{|x-a|^2}{2t}}\;dx,&&t>0,
\end{array}\right.$$
where $\delta_a$ denotes the point measure with respect to $a$.
In the univariate case we are able to show a general definability result with parameters.
The proof is rather easy and relies on the above mentioned result by Speissegger.

\vs{0.5cm}
{\bf Theorem A}

\vs{0.1cm}
{\it Let $\ma{M}$ be an o-minimal structure on the real field. Let $p\in \IN_0$ and let $A\subset \IR^p\times \IR$ be definable in $\ma{M}$. Then the parameterized integral 
	$$\IR^p\times \IR\times \IR_{>0}\to [0,1], (u,a,t)\mapsto P_a(B_t\in A_u),$$
	is definable in the Pfaffian closure $\ma{P}(\ma{M})$ of $\ma{M}$.}

\vs{0.5cm}
In the multivariate case we concentrate on the Brownian motion with start value zero (i.e. the standard Brownian motion) and work in the globally subanalytic setting.
Let $n\geq 2$ and let $A$ be a globally subanalytic subset of $\IR^n$.
Away from $0$ we obtain a definability result.

\vs{0.5cm}
{\bf Theorem B}

\vs{0.1cm}
{\it Let $a\in\IR_{>0}$. The restriction $[a,\infty[\;\to [0,1], t\mapsto P_0(B_t\in A),$  
is constructible. In particular, it is definable in the o-minimal expansion $\IR_{\an,\exp}$.}

\vs{0.5cm}
At $0$ we are able to show an  interesting asymptotic expansion.

\vs{0.5cm}
{\bf Theorem C}

\vs{0.1cm}
{\it There are $q\in \IN$ and for $k\in \IN_0$ and $l\in\{0,\ldots,n-2\}$ there are $a_{kl}\in \IR$ where $a_{01}=\ldots=a_{0,n-2}=0$ such that
$$P_0(B_t\in A)\sim \sum_{k=0}^\infty \big(a_{k0}+a_{k1}\log t+\ldots+a_{k,n-2}(\log t)^{n-2}\big)t^\frac{k}{q}.$$}

\vs{0.5cm}
The above type of series gives the right asymptotic scale for the time evolution $t\mapsto P_0(B_t\in A)$. Namely, the asymptotic expansion of the latter vanishes if and only if the given globally subanalytic set is thin at the origin.

Asymptotic
expansions of the above type often occur in analysis: for the context of the Riemann mapping theorem see Lehmann [19] and [12], for the context of the Dirichlet problem see Wasow [23] and [13].

In the situation of Theorem C the asymptotic expansion is in general not convergent. In particular, the function $t\mapsto P_0(B_t\in A)$ is in general not definable in $\IR_{\an,\exp}$. 
Although we have the same kind of asymptotic expansion as in the above mentioned situations of the Riemann mapping theorem or the Dirichlet problem, the time evolution of the probability is in contrast to them in general not in the quasianalytic Ilyashenko class (see Ilyashenko and Yakovenko [11, Section 24], [12, 13] and Speissegger [22] for definitions).
In general, logarithmic terms appear in the asymptotic expansion. But even if no logarithmic terms occur, the time evolution is in general not definable in the o-minimal structure $\IR_\ma{Q}$ (see [16]) related to the quasianalytic Ilyashenko class.

\section*{Notations}

By $\IN=\{1,2,3,\ldots\}$ we denote the set of natural numbers and by $\IN_0=\{0,1,2,3,\ldots\}$ the set of natural numbers with $0$.
We let $\IR_{>0}:=\;]0,\infty[\;=\{x\in\IR\mid x>0\}$ be the set of positive real
numbers. 
Given a function $f:X\to Y$ and a subset $U$ of $X$, we denote by $f|_U$ the restriction of $f$ to $U$.
Given a subset $C$ of $X\times Y$ and $x\in X$, we set $C_x:=\{y\in Y\mid (x,y)\in C\}$.

By $|x|$ we denote the euclidean absolute value of $x\in \IR^n$.
Given $a\in \IR^n$ and $r\in \IR_{>0}$ we set $B(a,r):=\{x\in \IR^n\mid |x-a|<r\}$, the euclidean open ball with center $a$ and radius $r$.
Let $A$ be a subset of $\IR^n$. The distance of $x\in\IR^n$ to $A$ is denoted by $\mathrm{dist}(x,A)$.

Let $f,g:\;]0,\infty[\;\to \IR$ be functions. We write $f=O(g)$ as $x\to 0$ if there are constants $a,c\in \IR_{>0}$ such that
$|f(x)|\leq c|g(x)|$ for all $x\in\;]0,a[$. We write $f=o(g)$ if for every $c\in\IR_{>0}$ there is $a\in\IR_{>0}$ such that
$|f(x)|\leq c|g(x)|$ for all $x\in \;]a,c[$. Likewise for $x\to \infty$. 

By $\IR[T_1,\ldots,T_k]$ we denote the polynomial ring in $k$ variables over the real field.

\vs{0.3cm}
Given $n\in \IN$ and a Borel subset $X$ of $\IR^n$ we set
$$\Phi_X:\IR_{>0}\to [0,1], t\mapsto \frac{1}{(2\pi t)^{\frac{n}{2}}}\int_Xe^{-\frac{|x|^2}{2t}}\;dx.$$

\section{Preliminaries}

We assume basic knowledge of o-minimality (see for example Van den Dries [7]) and of globally subanalytic sets and functions (see for example Van den Dries and Miller [10]).

\vs{0.5cm}
Cluckers and D. Miller [2] have introduced the notion of constructible functions:
Let $A$ be a globally subanalytic subset of some $\IR^n$. A function $f:A\to \IR$ is {\bf constructible} if there are $l\in \IN_0$, a polynomial $P\in \IR[X_1,\ldots,X_l,Y_1,\ldots,Y_l]$ and globally subanalytic functions
$\varphi_1,\ldots,\varphi_l:\IR^n\to \IR$ and $\psi_1,\ldots,\psi_l:\IR^n\to \IR_{>0}$ such that 
$$f=P\big(\varphi_1,\ldots,\varphi_l,\log(\psi_1),\ldots,\log(\psi_l)\big).$$

\vs{0.5cm}
{\bf 1.1 Remark}

\vs{0.1cm}
{\it Let $f:A\to \IR$ be constructible. Then $f$ is definable in $\IR_{\an,\exp}$.} 

\vs{0.5cm}
Let $p\in \IN_0$ be the total degree of the above polynomial $P$ with respect to the variable $Y=(Y_1,\ldots,Y_l)$.
We say that $f$ is {\bf constructible of logarithmic power at most $p$}.

\vs{0.5cm}
{\bf 1.2 Proposition}

\vs{0.1cm}
{\it Let $f:\;]0,\infty[\;\to \IR$ be constructible of logarithmic power at most $p$.
\begin{itemize}
	\item[(1)] There is a sufficiently small $a$ and there are $q\in \IN$ and Puiseux series $P_j(x)=\sum_{k=m_j}^\infty c_{jk}x^{k/q}$ (where $m_j\in \IZ$) for $j\in \{0,\ldots,p\}$ which converge on $]0,a[$ such that on $]0,a[$ 
	$$f(x)=\sum_{j=0}^pP_j(x)\big(\log x\big)^j.$$
	\item[(2)] There is a sufficiently large $b$ and there are $q\in \IN$ and Puiseux series $Q_j(x)=\sum_{k=m_j}^\infty c_{jk}x^{-k/q}$ (where $m_j\in\IZ$) for $j\in \{0,\ldots,p\}$ which converge on $]b,\infty[$ such that on $]b,\infty[$
	$$f(x)=\sum_{j=0}^pQ_j(x)\big(\log x\big)^j.$$
\end{itemize}}
{\bf Proof:}

\vs{0.1cm}
This follows by developing unary globally subanalytic functions into Puiseux series (see for example Van den Dries [6]) and by applying the properties of the logarithm. 
\hfill$\Box$

\vs{0.5cm}
From the work of Cluckers and D. Miller [2, 3, 4] we obtain the following:

\vs{0.5cm}
{\bf 1.3 Fact} 

\vs{0.1cm}
{\it Let $f:\IR^{m+n}\to \IR$ be constructible of logarithmic power at most $p$. 
	Assume that for every $x\in \IR^m$ the
	function 
	$\IR^n\to \IR, y\mapsto f(x,y),$
	is integrable.
	Then the function
	$$\IR^m\to \IR, x\mapsto \int_{\IR^n}f(x,y)\;dy,$$
	is constructible of logarithmic power at most $p+n$.}

\vs{0.5cm}
{\bf 1.4 Definition}

\begin{itemize}
	\item[(a)]  
	By $\ma{A}$ we denote the ring consisting of all formal series of the form $\sum_{k=0}^\infty g_k(\log x)x^{k/q}$
	where $q\in \IN$ and $g_k\in \IR[T]$ for all $k\in \IN_0$ with $g_0$ constant.\\
	Let $p\in \IN_0$. By $\ma{A}_p$ we denote the additive subgroup of $\ma{A}$ consisting of all $\sum_{k=0}^\infty g_k(\log x)x^{k/q}\in \ma{A}$
	such that $\deg(g_k)\leq p$ for all $k\in\IN$.
	By $\ma{A}_p^{\mathrm{conv}}$ we denote the subgroup of $\ma{A}_p$ consisting of all $F\in \ma{A}_p$ that converge for sufficiently small $x$.
	Moreover, we denote by $\ma{A}_{-1}=\IR$ the additive subgroup of $\ma{A}$ consisting of all constant series.
    \item[(b)]  
    By $\ma{B}$ we denote the ring consisting of all formal series of the form $\sum_{k=0}^\infty g_k(\log x)x^{-k/q}$
    where $q\in \IN$ and $g_k\in \IR[T]$ for all $k\in \IN_0$ with $g_0$ constant.\\
    Let $p\in \IN_0$. By $\ma{B}_p$ we denote the additive subgroup of $\ma{B}$ consisting of all $\sum_{k=0}^\infty g_k(\log x)x^{-k/q}\in \ma{B}$
    such that $\deg(g_k)\leq p$ for all $k\in\IN$.
    By $\ma{B}_p^{\mathrm{conv}}$ we denote the subgroup of $\ma{B}_p$ consisting of all $F\in \ma{B}_p$ that converge for sufficiently large $x$.
\end{itemize} 

\vs{0.2cm}
{\bf 1.5 Proposition}

\vs{0.1cm}
{\it Let $f:\;]0,\infty[\;\to \IR$ be a function and let $p\in\IN_0$.
\begin{itemize}
\item[(1)] The following are equivalent:
\begin{itemize}
\item[(i)] There is $a>0$ such that $f|_{]0,a[}$ is bounded and constructible of logarithmic power at most $p$.
\item[(ii)] There is $F\in \ma{A}_p^{\mathrm{conv}}$ such that $f(x)=F(x)$ for all sufficiently small $x$.
\end{itemize}
\item[(2)] The following are equivalent:
\begin{itemize}
\item[(i)] There is $b>0$ such that $f|_{]b,\infty[}$ is bounded and constructible of logarithmic power at most $p$.
\item[(ii)] There is $F\in \ma{B}_p^{\mathrm{conv}}$ such that $f(x)=F(x)$ for all sufficiently large $x$.
\end{itemize}
\end{itemize}}
{\bf Proof:}

\vs{0.1cm}
We show (1):

\vs{0.2cm} 
(i) $\Rightarrow$ (ii): Let $a>0$ such that $f|_{]0,a[}$ is bounded and constructible of logarithmic power at most $p$. By Proposition 1.2 we see after choosing a smaller $a$ if necessary that there are $q\in \IN$ and convergent Puiseux series $P_j(x)=\sum_{k=m_j}^\infty c_{jk}x^{k/q}$
for $j\in \{0,\ldots,p\}$ that converge on $]0,a[$ such that $f(x)=\sum_{j=0}^p P_j(x)\big(\log x\big)^j$ on $]0,a[$.
Since $f$ is bounded in a neighbourhood of $0$ we can assume that $m_0=0$ and $m_j>0$ for $j\in \{1,\ldots,p\}$.
Setting 
$$F(x):=\sum_{k=0}^\infty \big(c_{0k}+c_{1k}\log x+\ldots+c_{pk}(\log x)^p\big)x^{\frac{k}{q}}$$
we are done.

\vs{0.2cm}
(ii) $\Rightarrow$ (i): Let $a>0$ and let $F=\sum_{k=0}^\infty g_k(\log x)x^{k/p}\in \ma{A}_p^\mathrm{conv}$ such that $f(x)=F(x)$ for all $x\in ]0,2a[$. 
For $k\in \IN_0$ let $g_k(T)=d_{k0}+d_{k1}T+\ldots+d_{kp}T^p$.
For $j\in \{0,\ldots,p\}$ set $P_j(x):=\sum_{k=0}^\infty d_{kj}x^{k/q}$. Then $P_j(x)$ is a Puiseux series that converges on $]0,2a[$. Hence it defines a globally subanalytic function $\varphi_j:\;]0,2a[\;\to \IR$.
Let 
$$P:=X_0+X_1Y+\ldots+X_pY^p\in \IR[X_0,\ldots,X_p,Y].$$
Then $f|_{]0,a[}=P(\varphi_0,\ldots,\varphi_p,\log x)$ and therefore $f$ is constructible on $]0,a[$. That $f$ is bounded on $]0,a[$ follows from the fact that each $P_j$ respectively $\varphi_j$ is bounded on $]0,a[$ and that
$d_{0j}=0$ for $j\in \{1,\ldots,p\}$.
\hfill$\Box$

\vs{0.5cm}
{\bf 1.6 Definition}

\vs{0.1cm}
Let $f:\;]0,\infty[\;\to \IR$ be a function.
\begin{itemize}
	\item[(a)] We say that $f$ has {\bf asymptotic expansion} $F=\sum_{k=0}^\infty g_k(\log x)x^{k/q}\in \ma{A}$ at $0$ and write $f\sim F$ as $x\to 0$
	if for every $N\in \IN_0$ 
	$$f(x)-\sum_{k=0}^N g_k(\log x)x^{\frac{k}{q}}=o(x^{\frac{N}{q}})$$
	as $x\to 0$.
	\item[(b)] We say that $f$ has {\bf asymptotic expansion} $F=\sum_{k=0}^\infty g_k(\log x)x^{-k/q}\in \ma{B}$ at $\infty$ and write $f\sim F$ as $x\to \infty$
	if for every $N\in \IN_0$ 
	$$f(x)-\sum_{k=0}^N g_k(\log x)x^{-\frac{k}{q}}=o(x^{-\frac{N}{q}})$$
	as $x\to \infty$.
\end{itemize}

\vs{0.2cm}
{\bf 1.7 Remark}

\vs{0.1cm}
{\it The asymptotic expansion at $0$ respectively at $\infty$ is uniquely determined if it exists.}

\vs{0.5cm}
{\bf 1.8 Remark}

\vs{0.1cm}
{\it Let $f:\IR_{>0}\to \IR$ be bounded and constructible of logarithmic power at most $p$.
\begin{itemize}
	\item[(1)] There is $F\in \ma{A}_p^\mathrm{conv}$ such that $f\sim F$ at $0$.
	\item[(2)] There is $F\in \ma{B}_p^\mathrm{conv}$ such that $f\sim F$ at $\infty$.
\end{itemize}}
{\bf Proof:}

\vs{0.1cm}
This follows from Proposition 1.5.
\hfill$\Box$

\section{The univariate case}

In the one-dimensional case we obtain a definability result including parameters.
Let $\mathcal{M}$ be an o-minimal expansion of the real field. We show Theorem A of the introduction.
By $\ma{P}(\ma{M})$ we denote the Pfaffian closure of $\ma{M}$; note that $\ma{P}(\ma{M})$ is again o-minimal (see Speissegger [21]).

\vs{0.5cm}
{\bf 2.1 Theorem}

\vs{0.1cm}
{\it Let $\mathcal{M}$ be an o-minimal structure on the real field. Let $p\in\IN_0$ and let $A\subset \IR^p\times \IR$ be definable in $\mathcal{M}$.
Then the parameterized integral
\begin{eqnarray*}
\chi=\chi_A: \IR^p \times \IR\times \IR_{>0} & \longrightarrow & [0,1],\\
(u,a,t) & \longmapsto & P_a(B_t\in A_u)=
\frac{1}{\sqrt{2 \pi t}}\int\limits_{A_u} e^{\frac{-(x-a)^2}{2t}}\; dx.
\end{eqnarray*} 
is definable in the Pfaffian closure $\mathcal{P}(\mathcal{M})$ of $\mathcal{M}$.}

\vs{0.1cm}
{\bf Proof:}

\vs{0.1cm}
We start the proof with the following

\vs{0.1cm}
{\bf Claim:}
The error function
$$\mathrm{erf}:\IR\to \IR, x\mapsto \frac{2}{\sqrt{\pi}}\int_0^x e^{-s^2}\;ds,$$
is definable in $\ma{P}(\ma{M})$.

\vs{0.1cm}
{\bf Proof of the claim:}
The exponential function is definable in $\ma{P}(\ma{M})$, hence $s\mapsto e^{-s^2}$ is definable in $\ma{P}(\ma{M})$.
Up to a multiplicative factor, the error function is the antiderivative of the latter. Hence we obtain the claim by [21, p. 190].
\hfill$\Box_{\mathrm{Claim}}$

\vs{0.2cm}
By additivity of the integral and by cell decomposition, we may assume that $A$ is a cell. Moreover, we may assume that $A$ is of type band over its basis $\pi(A)$ where $\pi:\IR^p\times\IR\to \IR^p, (u,x)\mapsto u,$ denotes the projection onto the first factor.
Hence there are continuous definable functions $f,g:\pi(A)\to \IR\cup\{\pm \infty\}$ with $f<g$ such that
$$A_u=\left\{\begin{array}{ccc}
\big\{x\in\IR\mid f(u)<x<g(u)\big\},&&a\in \pi(A),\\
&\mbox{if}&\\
\emptyset,&&u\notin \pi(A).\\
\end{array}\right.$$
By the substitution $y=(x-a)/\sqrt{2t}$ and the main theorem of calculus we obtain 
that
$$\chi_A(u,a,t)=\left\{\begin{array}{ccc}
\frac{1}{2}\left(
\mathrm{erf}\left(\frac{g(u)-a}{\sqrt{2t}} \right)-\mathrm{erf} \left(\frac{f(u)-a}{\sqrt{2t}}\right)\right),&&u\in\pi(A), \\
&\mbox{if}&\\
0,&&a\notin \pi(A).\\
\end{array}\right.$$
(Note that $\mathrm{erf}(+\infty)=1$ and $\mathrm{erf}(-\infty)=-1$.) Hence we see by the claim that $\chi_A$ is definable in the Pfaffian closure $\mathcal{P}(\mathcal{M})$ of $\mathcal{M}$.
\hfill$\Box$

\vs{0.5cm}
We focus on the start value zero and consider the case of a single definable subset of $\IR$, i.e. a finite union of intervals and points. 
The results below will be generalized in Section 3 to the multivariate case (compare with Theorems B and C of the introduction).

\vs{0.5cm}
{\bf 2.2 Proposition}

\vs{0.1cm}
{\it Let $A$ be a finite union of intervals and points. The following holds:
\begin{itemize}
	\item[(1)] For every $a\in \IR_{>0}$ the restriction $\Phi_A|_{[a,\infty[}$ is globally subanalytic.
	\item[(2)] There is $F\in \ma{A}_{-1}$ such that $\Phi_A(t)\sim F$ as $t\to 0$. We have that $F\neq 0$ if and only if $0\in \overline{\mathring{A}}$.
\end{itemize}}
{\bf Proof:}

\vs{0.1cm}
Since $A\setminus \mathring{A}$ is a Lebesgue nullset we may assume that $A$ is open.
Hence $A$ is a finite union of disjoint open interval.

\vs{0.2cm}
(1): By the linearity of the integral we may assume that $A$ is a nonempty open interval of the form $]c,d[$ where
$-\infty\leq c<d\leq\infty$.
We obtain
$$\Phi_A(t)=\frac{1}{2}\left(
\mathrm{erf}\left(\frac{d}{\sqrt{2t}} \right)-\mathrm{erf} \left(\frac{c}{\sqrt{2t}}\right)\right).$$
If $d=\infty$ then the first summand is constant $1$. If $d<\infty$ we may assume for simplicity that $d\geq 0$.
The restriction of the error function to the interval $[0,d/\sqrt{2a}]$ is restricted analytic. Hence the first summand restricted to $[a,\infty[$ is globally subanalytic. The second summand is treated similarly.

\vs{0.2cm}
(2): Again by the linearity if the integral it is enough to show the following:
Let $A$ be a nonempty open interval. If $0\in\overline{A}$ there is $F\in \IR\setminus\{0\}$ such that $\Phi_A\sim F$ as $t\to 0$. If $0\notin\overline{A}$ then $\Phi_A\sim 0$ as $t\to 0$.
 
To show the first statement it is enough by the positivity of the integrand to assume that $A$ is of the form $]0,b[$ for some $b>0$.
We have
$$\Phi_A(t)=\frac{1}{2}
\mathrm{erf}\Big(\frac{b}{\sqrt{2t}}\Big).$$
Since $\mathrm{erf}(x)=1+o(x^{-N})$ as $x\to \infty$ for every $N\in\IN$ (see for example Van den Dries at al. [9, (5.7)]) we get that $\Phi_A(t)\sim 1/2$.
To show the second statement it is enough to assume that $A$ is of the form $]b,\infty[$ for some $b>0$.
We have 
$$\Phi_A(t)=\frac{1}{2}\Big(1-
\mathrm{erf}\Big(\frac{b}{\sqrt{2t}}\Big)\Big).$$
Again using that $\mathrm{erf}(x)=1+o(x^{-N})$ as $x\to \infty$ for every $N\in\IN$ we obtain that
$\Phi_A(t)\sim 0$.
\hfill$\Box$

\section{The multivariate case}

From now on we assume that $n\geq 2$ and assume that $A$ is a globally subanalytic subset of $\IR^n$. We concentrate on the Brownian motion with start value $0$.

We introduce some notation we will use throughout this section.

\vs{0.5cm}
Let 
\begin{eqnarray*}
	\Theta_n:\IR_{>0}\times \;]-\pi,\pi[\;\times \;]-\pi/2,\pi/2[^{n-2}&\longrightarrow& \IR^n,\\
	(r,\varphi,\vartheta)&\longmapsto& \Theta_n(r,\varphi,\vartheta),
\end{eqnarray*}
be the $n$-dimensional polar coordinates where $\vartheta:=(\vartheta_1,\ldots,\vartheta_{n-2})$.
Since the map $\Theta_n$ is globally subanalytic we have that $B:=\Theta^{-1}(A)$ is also globally subanalytic.
By the transformation formula we obtain
$$\Phi_A(t)=\frac{1}{(2\pi t)^{\frac{n}{2}}}\int_{B}e^{-\frac{r^2}{2t}}r^{n-1}C_n(\vartheta)\,d(r,\varphi,\vartheta)$$
where 
$$C_n:\;]-\pi/2,\pi/2[^{n-2}\to \IR_{>0}, \vartheta\mapsto \prod_{j=1}^{n-2}(\cos \vartheta_j)^j.$$
For $r>0$ let 
$$B_r:=\Big\{(\varphi,\vartheta)\in \;]-\pi,\pi[\;\times\; ]-\pi/2,\pi/2[^{n-2}\;\Big\vert\;
(r,\varphi,\vartheta)\in B\Big\}.$$
By Fubini's Theorem, we obtain
\begin{eqnarray*}
	\Phi_A(t)&=&
	\frac{1}{(2\pi t)^{\frac{n}{2}}}\int_{\IR_{>0}}e^{-\frac{r^2}{2t}} r^{n-1}\Big(\int_{B_r}
	C_n(\vartheta)\;d(\varphi,\vartheta)\Big)\, dr\\
	&=&\frac{1}{(2\pi t)^{\frac{n}{2}}}\int_{\IR_{>0}}e^{-\frac{r^2}{2t}} r^{n-1}\Delta(r)\, dr
\end{eqnarray*}
where
$$\Delta:\IR_{>0}\to \IR_{\geq 0}, r\mapsto \int_{B_r}C_n(\vartheta)\;d(\varphi,\vartheta).$$
Note that $\Delta$ is bounded.
Integrating first with respect to $\varphi$ does not give a logarithmic term.
By Fact 1.3  we obtain that $\Delta$ is constructible of logarithmic power at most $n-2$.

\subsection{Results for time approaching infinity}

We establish a more detailed version of Theorem B.

\vs{0.5cm}
{\bf 3.1 Theorem}

\vs{0.1cm}
{\it Let $a\in \IR_{>0}$. Then $\Phi_A|_{[a,\infty[}$ is constructible of logarithmic power at most $n-1$.}

\vs{0.1cm}
{\bf Proof:}

\vs{0.1cm}
By Proposition 1.5(2) there are 
$\alpha>0$ and $F=\sum_{k=0}^\infty g_k(\log r)r^{-k/q}\in \ma{B}_{n-2}^\mathrm{conv}$ such that $\Delta(r)=F(r)$
on $[\alpha/2,\infty[$.
By the linearity of the integral the statement of the theorem follows from the following two claims:

\vs{0.2cm}
{\bf Claim 1:}
The function  
$$\Sigma_1:[a,\infty[\;\to \IR, t\mapsto\frac{1}{(2\pi t)^{\frac{n}{2}}}\int_\alpha^\infty e^{-\frac{r^2}{2t}}r^{n-1}\Delta(r)\;dr,$$
is constructible of logarithmic power at most $n-1$.

\vs{0.2cm}
{\bf Claim 2:}
The function
$$\Sigma_2:[a,\infty[\;\to \IR,t\mapsto\frac{1}{(2\pi t)^{\frac{n}{2}}}\int_0^\alpha e^{-\frac{r^2}{2t}}r^{n-1}\Delta(r)\;dr,$$
is constructible of logarithmic power at most $n-1$.

\vs{0.2cm}
{\bf Proof of Claim 1:}
By the additivity of the integral we may assume that $\Delta$ respectively $F$ is of the form $\Delta(r)=Q(r)(\log r)^l$ where $l\in \{0,\ldots,n-2\}$ and $Q(r)=\sum_{k=0}^\infty c_kr^{-k/q}$ is a Puiseux series with $c_0=0$ if $l>0$ which converges on $[\alpha/2,\infty[$.
We obtain by substituting $s=r/\sqrt{2t}$
\begin{eqnarray*}
	\Sigma_1(t)&=&\frac{1}{(2\pi t)^{\frac{n}{2}}}\int_\alpha^\infty e^{-\frac{r^2}{2t}}r^{n-1}Q(r)\big(\log r\big)^l\;dr\\
	&=&\pi^{-\frac{n}{2}}\int_{\alpha/\sqrt{2t}}^\infty e^{-s^2}s^{n-1}Q(s\sqrt{2t})\big(\log s\sqrt{2t}\big)^l\;ds.\\
\end{eqnarray*}
We set
$$\Sigma_{11}:[a,\infty[\;\to \IR, t\mapsto \pi^{-\frac{n}{2}}\int_{\alpha/\sqrt{2t}}^{\alpha/\sqrt{2a}} e^{-s^2}s^{n-1}Q(s\sqrt{2t})\big(\log s\sqrt{2t}\big)^l\;ds$$
and 
$$\Sigma_{12}:[a,\infty[\;\to \IR, t\mapsto \pi^{-\frac{n}{2}}\int_{\alpha/\sqrt{2a}}^\infty e^{-s^2}s^{n-1}Q(s\sqrt{2t})\big(\log s\sqrt{2t}\big)^l\;ds.$$
The integrand
$$[a,\infty[\;\times\; [0,\alpha/\sqrt{2a}]\to \IR, (t,s)\mapsto e^{-s^2}s^{n-1}Q(s\sqrt{2t})\big(\log s\sqrt{2t}\big)^l,$$
of $\Sigma_{11}$ is constructible of logarithmic power at most $n-2$. By Fact 1.3 we get that $\Sigma_{11}$ is constructible of logarithmic power at most $n-1$.\\
The integrand of the second function is not constructible. We compute $\Sigma_{12}$ directly:
\begin{eqnarray*}
\Sigma_{12}(t)&=&\pi^{-\frac{n}{2}}\int_{\alpha/\sqrt{2a}}^\infty e^{-s^2}s^{n-1}Q(s\sqrt{2t})\big(\log s\sqrt{2t}\big)^l\;ds\\
&=&\pi^{-\frac{n}{2}}\int_{\alpha/\sqrt{2a}}^\infty e^{-s^2}s^{n-1}\Big(\sum_{k=0}^\infty c_k(s\sqrt{2t})^{-\frac{k}{q}}\Big)\big(\log s\sqrt{2t}\big)^l\;ds\\
&=&\pi^{-\frac{n}{2}}\int_{\alpha/\sqrt{2a}}^\infty e^{-s^2}s^{n-1}\Big(\sum_{k=0}^\infty c_k (s\sqrt{2t})^{-\frac{k}{q}}\Big(\sum_{p=0}^l\frac{\binom{l}{p}}{2^p}\big(\log\sqrt{2}s\big)^{l-p}(\log t)^p\Big)\Big)\;ds\\
&=&\sum_{k=0}^\infty\sum_{p=0}^lC_{kp}\Big(\int_{\alpha/\sqrt{2a}}^\infty e^{-s^2}s^{-\frac{k}{q}+n-1}\big(\log\sqrt{2}s\big)^{l-p}\;ds\Big)t^{-\frac{k}{2q}}(\log t)^p
\end{eqnarray*}
where $$C_{kp}:=c_k\frac{\binom{l}{p}}{2^{p+\frac{k}{2q}}\pi^{\frac{n}{2}}}.\;\;\;\;\;\;\;\;\;\;\;\;\;\;\;\;\;\;\;\;\;(*)$$
For $(k,p)\in
\IN_0\times\{0,\ldots,l\}$ set
$$D_{kp}:= \int_{\alpha/\sqrt{2a}}^\infty e^{-s^2}s^{-\frac{k}{q}+n-1}\big(\log\sqrt{2}s\big)^{l-p}\;ds.\;\;\;\;\;\;\;\;\;\;\;\;\;\;\;\;\;\;\;\;\;(**)$$
There is $c>0$ such that $\big\vert\big(\log\sqrt{2}s\big)^j\big\vert\leq cs$ for all $s\geq \alpha/\sqrt{2a}$ and all $j\in \{0,\ldots,n-2\}$.
We obtain
\begin{eqnarray*}
|D_{kp}|&\leq& c\int_{\alpha/\sqrt{2a}}^\infty e^{-s^2}s^{-\frac{k}{q}+n}\;ds\\
&\leq&c\Big(\frac{\alpha}{\sqrt{2a}}\Big)^{-\frac{k}{q}}\int_{\alpha/\sqrt{2a}}^\infty e^{-s^2}s^n\;ds.
\end{eqnarray*}
Hence
$$|D_{kp}|\leq cd\Big(\frac{\alpha}{\sqrt{2a}}\Big)^{-\frac{k}{q}}$$
for all $(k,p)\in \IN_0\times \{0,\ldots,l\}$ where
$$d:=\int_{\alpha/\sqrt{2a}}^\infty e^{-s^2}s^n\;ds\in \IR_{>0}.$$
Since the Puiseux series $P(r)=\sum_{k=0}^\infty c_kr^{-k/q}$ is convergent on $[\alpha/2,\infty[$ we find some $C>0$ such that 
$$|c_k|\leq C\big(\frac{\alpha}{2}\big)^{\frac{k}{q}}.\;\;\;\;\;\;\;\;\;\;\;\;\;\;\;\;\;\;\;\;\;(***)$$
for all $k\in \IN_0$.
By $(*) - (***)$ we see that
$$\Sigma_{12}(t)=\sum_{k=0}^\infty\sum_{p=0}^lC_{kp}D_{kp}t^{-\frac{k}{2q}}(\log t)^p\in \ma{B}_{n-2}^\mathrm{conv}$$
converges on $[a,\infty[$. We obtain by Proposition 1.5(2) that $\Sigma_{12}$ is constructible of logarithmic power at most $n-2$.\\
Since $\Sigma_1=\Sigma_{11}+\Sigma_{12}$ we get that $\Sigma_1$ is constructible of logarithmic power at most $n-1$.
\hfill$\Box_{\mathrm{Claim\;1}}$

\vs{0.2cm}
{\bf Proof of Claim 2:}
Since the function 
$$[a,\infty[\;\times\; [0,\alpha]\to \IR,(t,r)\mapsto e^{-\frac{r^2}{2t}},$$
is globally subanalytic we get that the integrand
$$[a,\infty[\;\times\;[0,\alpha]\to \IR, (t,r)\mapsto  e^{-\frac{r^2}{2t}}r^{n-1}\Delta(r),$$
of $\Sigma_2$ is constructible of logarithmic power at most $n-2$. By Fact 1.3 we get that $\Sigma_2$ is constructible of logarithmic power at most $n-1$.
\hfill$\Box_{\mathrm{Claim\;2}}$

\vs{0.2cm}
Since $\Phi_A|_{[a,\infty[}=\Sigma_1+\Sigma_2$ we are done by Claims 1 and 2.
\hfill$\Box$

\vs{0.5cm}
Note that Theorem 3.1 resp. Theorem B does not follow from the work of Comte, Lion, Rolin [5, 20] and of Cluckers, D. Miller [2, 3, 4] since the integrand of $\Phi_A$ is for unbounded $A$ not globally subanalytic or constructible.

\vs{0.5cm}
{\bf 3.2 Corollary}

\vs{0.1cm}
{\it Let $a\in\IR_{>0}$. Then $\Phi_A|_{[a,\infty[}$ is definable in $\IR_{\an,\exp}$.}

\vs{0.1cm}
{\bf Proof:}

\vs{0.1cm}
This follows from Theorem 3.1 and Remark 1.1.
\hfill$\Box$

\vs{0.5cm}
{\bf 3.3 Corollary}

\vs{0.1cm}
{\it There is $F\in \ma{B}_{n-1}^{\mathrm{conv}}$ such that
	$\Phi_A\sim F$ at $\infty$.}

\vs{0.1cm}
{\bf Proof:}

\vs{0.1cm}
By Theorem 3.1 we know that $\Phi_A|_{[1,\infty[}$ is constructible of logarithmic power at most $n-1$.
Since $\Phi_A(t)\in [0,1]$ for all $t\in \IR_{>0}$ we have that $\Phi_A$ is bounded. Remark 1.8(2) gives the result.
\hfill$\Box$

\vs{0.5cm}
The next example shows that logarithmic terms may occur.

\vs{0.5cm}
{\bf 3.4 Example}

\vs{0.1cm}
{\it Consider the globally subanalytic subset $A$ of $\IR^2$ which is given in polar coordinates
	by
	$$\Big\{(r,\varphi)\in \IR_{>0}\times\; ]-\pi,\pi[\;\Big\vert\; r>1, 0<\varphi<1/r^2\Big\}.$$
	Then the asymptotic expansion of $\Phi_A$ at $\infty$ contains logarithmic terms.}

\vs{0.1cm}
{\bf Proof:}
We obtain for $r\in \;]1,\infty[$ that
$\Delta(r)=\int_{0}^{1/r}\;d\varphi=1/r^2$.
For $t\geq 1$ we get
\begin{eqnarray*}
	\Phi_A(t)
	&=&\frac{1}{2\pi t}\int_1^{\infty}e^{-\frac{r^2}{2t}}r^{-1}\;dr\\
	&=&\frac{1}{2\pi t}\int_{1/\sqrt{2t}}^\infty e^{-s^2}s^{-1}\;ds\\
	&=&\frac{1}{2\pi t}\Big(\int_{1/\sqrt{2t}}^1 e^{-s^2}s^{-1}\;ds+c\Big)\\
\end{eqnarray*}
where $c=\int_1^\infty e^{-s^2}s^{-1}\;ds$. We have
\begin{eqnarray*}
	\int_{1/\sqrt{2t}}^1 e^{-s^2}s^{-1}\;ds&=&\int_{1/\sqrt{2t}}^1\sum_{k=0}^\infty\Big(\frac{(-1)^k}{k!}s^{2k-1}\Big)\;ds\\
	&=&\sum_{k=0}^\infty\frac{(-1)^k}{k!}\int_{1/\sqrt{2t}}^1s^{2k-1}\;ds\\
	&=&\Big[\log s\Big]_{1/\sqrt{2t}}^1+\sum_{k=1}^\infty \frac{(-1)^k}{k!2k}\Big[s^{2k}\Big]_{1/\sqrt{2t}}^1\\
	&=&\frac{1}{2}\log t+\sum_{k=0}^\infty a_k t^{-k}\\
\end{eqnarray*}
where
$$a_0=\frac{1}{2}\log 2+\sum_{k=1}^\infty \frac{(-1)^k}{k!2k}$$
and $a_k=(-1)^{k-1}/(2^{k+1}k!k)$ for $k\geq 1$. We conclude that for $t\geq 1$
$$\Phi_A(t)=\frac{1}{2\pi}\Big(\frac{1}{2}\big(\log t\big)t^{-1}+\sum_{k=0}^\infty c_kt^{-k-1}\Big)$$
where $c_0=a_0+c$ and $c_k=a_k$ for $k\geq 1$.
Hence $\Phi_A|_{[1,\infty[}$ is constructible of logarithmic power exact $1$.
\hfill$\Box$

\subsection{Results for time approaching zero}

We consider now the case that time approaches $0$. We establish Theorem C.

\vs{0.5cm}
{\bf 3.5 Theorem}

\vs{0.1cm}
{\it There is $F\in \ma{A}_{n-2}$ such that $\Phi_A\sim F$ at $0$.}

\vs{0.1cm}
{\bf Proof:}

\vs{0.1cm}
By Proposition 1.5(1) there are $\alpha>0$ and $G=\sum_{k=0}^\infty g_k(\log r)r^{k/q}\in\ma{A}_{n-2}^{\mathrm{conv}}$ such that $\Delta(r)=G(r)$ on $]0,2\alpha[$. 
By the linearity of the integral 
the theorem follows by the following two claims:

\vs{0.2cm}
{\bf Claim 1:} Let 
$$\Pi_1:\;]0,\infty[\;\to \IR, t\mapsto \frac{1}{(2\pi t)^{\frac{n}{2}}}\int_0^\alpha e^{-\frac{r^2}{2t}}r^{n-1}\Delta(r)\;dr.$$
Then there is $F\in \ma{A}_{n-2}$ such that $\Pi_1\sim F$ at $0$.

\vs{0.2cm}
{\bf Claim 2:} Let
$$\Pi_2:\;]0,\infty[\;\to \IR, t\mapsto \frac{1}{(2\pi t)^{\frac{n}{2}}}\int_\alpha^\infty e^{-\frac{r^2}{2t}}r^{n-1}\Delta(r)\;dr.$$
For every $N\in \IN$ we have $\Pi_2(t)=o(t^N)$
as $t\to 0$.

\vs{0.2cm}
{\bf Proof of Claim 1:}
By the additivity of the integral we may assume that $\Delta$ respectively $G$ is of the form $\Delta(r)=P(r)(\log r)^l$ where $l\in \{0,\ldots,n-2\}$ and $P(r)=\sum_{k=0}^\infty c_kr^{k/q}$ is a Puiseux series with $c_0=0$ if $l>0$  which converges on $[0,2\alpha[$.
Fix $K\in \IN$. Then there are $\sigma>0$ with $K/q<\sigma<(K+1)/q$ and $C>0$ such that $|\sum_{k=K+1}^\infty c_kr^{k/q}(\log r)^l|\leq Cr^\sigma$ for all $r\in [0,\alpha]$.
Let 
$$\Pi_{1K1}:\;]0,\infty[\;\to \IR, t\mapsto \frac{1}{(2\pi t)^{\frac{n}{2}}}\int_0^\alpha e^{-\frac{r^2}{2t}}r^{n-1}\Big(\sum_{k=0}^K c_kr^{\frac{k}{q}}\Big)(\log r)^l\;dr,$$
and 
$$\Pi_{1K2}:\;]0,\infty[\;\to \IR, t\mapsto\frac{1}{(2\pi t)^{\frac{n}{2}}}\int_0^\alpha e^{-\frac{r^2}{2t}}r^{n-1}\Big(\sum_{k=K+1}^\infty c_kr^{\frac{k}{q}}\Big)(\log r)^l\;dr.$$
We obtain by substituting $s=r/\sqrt{2t}$
\begin{eqnarray*}
	\Pi_{1K1}(t)&=&\frac{1}{(2\pi t)^{\frac{n}{2}}}\int_0^\alpha e^{-\frac{r^2}{2t}}r^{n-1}\Big(\sum_{k=0}^K c_kr^{\frac{k}{q}}\Big)(\log r)^l\;dr\\
&=&\pi^{-\frac{n}{2}}\int_0^{\alpha/\sqrt{2t}} e^{-s^2}s^{n-1}\Big(\sum_{k=0}^K c_k(s\sqrt{2t})^{\frac{k}{q}}\Big)\big(\log s\sqrt{2t}\big)^l\;ds\\
		&=&\pi^{-\frac{n}{2}}\int_0^{\alpha/\sqrt{2t}} e^{-s^2}s^{n-1}\Big(\sum_{k=0}^K c_k (s\sqrt{2t})^{\frac{k}{q}}\Big(\sum_{p=0}^l\binom{l}{p}\big(\log\sqrt{2}s\big)^{l-p}(\log t)^p\Big)\Big)\;ds\\
		&=&\sum_{k=0}^K\sum_{p=0}^lC_{kp}\Big(\int_0^{\alpha/\sqrt{2t}} e^{-s^2}s^{\frac{k}{q}+n-1}\big(\log\sqrt{2}s\big)^{l-p}\;ds\Big)t^{\frac{k}{2q}}(\log t)^p
	\end{eqnarray*}
where $$C_{kp}:=c_k\frac{2^{\frac{k}{2q}-p}}{\pi^{\frac{n}{2}}}\binom{l}{p}.$$
For $(k,p)\in
\{0,\ldots,K\}\times\{0,\ldots,l\}$ set
$$\sigma_{kp}:\IR_{>0}\to \IR, t\mapsto \int_0^{\alpha/\sqrt{2t}}e^{-s^2}s^{\frac{k}{q}+n-1}\big(\log\sqrt{2}s\big)^{l-p}.$$
We have that
$$\lim_{t\to 0}\sigma_{kp}(t)=\int_0^\infty e^{-s^2}s^{\frac{k}{q}+n-1}\big(\log\sqrt{2}s\big)^{l-p}=:D_{kp}$$ 
exists in $\IR$.

\vs{0.2cm}
{\bf Claim 1.1:}
For every $N\in\IN$ we have that $\sigma_{kp}(t)=D_{kp}+o(t^N)$ as $t\to 0$.

\vs{0.1cm}
{\bf Proof of Claim 1.1:}
Let 
$$\tau_{kp}:\;]0,\infty[\;\to \IR, t\mapsto \int_{\alpha/\sqrt{2t}}^\infty e^{-s^2}s^{\frac{k}{q}+n-1}\big(\log\sqrt{2}s\big)^{l-p}.$$
Let $N\in \IN$. We have to show that $\tau_{kp}=o(t^N)$ as $t\to 0$.
Choosing $t_0$ small enough we have that 
$$|e^{-s^2}s^{\frac{k}{q}+n-1}\big(\log\sqrt{2}s\big)^{l-p}|\leq e^{-s}$$
for all $s\geq \alpha/\sqrt{2t_0}$. Hence 
$$|\tau_{kp}(t)|\leq \int_{\alpha/\sqrt{2t}}^\infty e^{-s}\;ds=e^{-\frac{\alpha}{\sqrt{2t}}}$$
for all $t\leq t_0$. This gives the assertion.
\hfill$\Box_{\mathrm{Claim\;1.1}}$

\vs{0.2cm}
For $(k,p)\in \{0,\ldots,K\}\times \{0,\ldots,l\}$ let $E_{kp}:=C_{kp}D_{kp}$.
Using the above presentation of $\Pi_{1K1}$ as a finite sum we obtain by Claim 1.1 that 
$$\Pi_{1K1}(t)\sim \sum_{k=0}^K\sum_{p=0}^l E_{kp} t^{\frac{k}{2q}}(\log t)^p\;\;\;\;\;\;\;\;(*)$$
as $t\to 0$.\\
We estimate $\Pi_{1K2}$. Again by substituting $s=r/\sqrt{2t}$, we obtain
\begin{eqnarray*}
	\big\vert\Pi_{1K2}(t)\big\vert&=&\Big\vert\frac{1}{(2\pi t)^{\frac{n}{2}}}\int_0^\alpha e^{-\frac{r^2}{2t}}r^{n-1}\Big(\sum_{k=K+1}^\infty c_kr^{\frac{k}{q}}\Big)(\log r)^l\;dr\Big\vert\\
	&\leq &\frac{1}{(2\pi t)^{\frac{n}{2}}}\int_0^\alpha e^{-\frac{r^2}{2t}}r^{n-1}\Big\vert\Big(\sum_{k=K+1}^\infty c_kr^{\frac{k}{q}}\Big)(\log r)^l\Big\vert\;dr\\
	&\leq&\frac{C}{(2\pi t)^{\frac{n}{2}}}\int_0^\alpha e^{-\frac{r^2}{2t}}r^{\sigma+n-1}\;dr\\
	&=&\frac{C}{(2\pi t)^{\frac{n}{2}}}\int_0^{\alpha/\sqrt{2t}}e^{-s^2}(s\sqrt{2t})^{\sigma+n-1}\sqrt{2t}\;ds\\
	&=&D\Big(\int_0^{\alpha/\sqrt{2t}}e^{-s^2}s^{\sigma+n-1}\;ds\Big)t^{\frac{\sigma}{2}}\\
	&\leq& Et^{\frac{\sigma}{2}}
\end{eqnarray*}
where 
$D:=C2^{\frac{\sigma}{2}}/\pi^{\frac{n}{2}}$
and 
$$E:=D\Big(\int_0^\infty e^{-s^2}s^{\sigma+n-1}\;ds\Big).$$
This shows that 
$$\Pi_{1K2}(t)=O(t^{\frac{\sigma}{2}})\;\;\;\;\;\;(**)$$
as $t\to 0$. 
Since $\Pi_1=\Pi_{1K1}+\Pi_{1K2}$ and since $\sigma>K/q$ we obtain by  $(*)$ and $(**)$ that 
$$\Pi_1(t)=\sum_{k=0}^K\sum_{p=0}^l E_{kp} t^{\frac{k}{2q}}(\log t)^p+o(t^{\frac{K}{2q}})$$
as $t\to 0$.
Since $K$ was arbitrary we obtain Claim 1.
\hfill$\Box_{\mathrm{Claim\;1}}$

\vs{0.2cm}
{\bf Proof of Claim 2:}
Since $\alpha>0$ and since $\Delta$ is bounded we find a constant $C>0$ such that 
$$\big\vert e^{-\frac{r^2}{2t}}r^{n-1}\Delta(r)\big\vert\leq Ce^{-\frac{r}{2t}} $$
for all $r\geq \alpha$ and all sufficiently small $t$.
We obtain for those $t$ that
$$\big\vert\Pi_2(t)\big\vert\leq
\frac{1}{(2\pi t)^{\frac{n}{2}}}\int_\alpha^\infty e^{-\frac{r}{2t}}\;dr=
\frac{2t}{(2\pi t)^{\frac{n}{2}}}e^{-\frac{\alpha}{2t}}.$$
Hence we every $N\in \IN$ we obtain that $\Pi_2(t)=o(t^N)$ as $t\to 0$.
\hfill$\Box_{\mathrm{Claim\;2}}$

\vs{0.2cm}
Since $\Phi_A=\Pi_1+\Pi_2$ the theorem is proven by Claims 1 and 2.
\hfill$\Box$

\vs{0.5cm}
The nest result gives a geometric criterion for the non-vanishing of the asymptotic expansion of $\Phi_A(t)$ at $0$.
It particular shows that we are working a priori with a right asymptotic scale.

\vs{0.5cm}
{\bf 3.6 Theorem}

\vs{0.1cm}
{\it The following are equivalent:
\begin{itemize}
	\item[(i)] The asymptotic expansion of $\Phi_A(t)$ at $0$ does not vanish.
	\item[(ii)] The origin is in the closure of the interior of $A$. 
\end{itemize}}
{\bf Proof:}

\vs{0.1cm}
We use the notation of the proof of Theorem 3.5.

\vs{0.2cm}
(i) $\Rightarrow$ (ii):
Assume that $0\notin \overline{\mathring{A}}$.
Then we have that $\dim(B_r)<n$ for all sufficiently small $r$. This gives that $\Delta(r)=0$ for all sufficiently small $r$.
After choosing an appropriate $\alpha$ we get that $\Phi_A(t)=\Pi_2(t)$.
We are done by Claim 2 of the previous proof.

\vs{0.2cm}
(ii) $\Rightarrow$ (i):
Assume that $0\in\overline{\mathring{A}}$.
Then $\dim(B_r)=n$ for all sufficiently small $r$. Since $C_n(\vartheta)$ is continuous and positive we see that $\Delta(r)>0$ for all sufficiently small $r$.
Since $\Delta$ is constructible we find some small $\varepsilon>0$ and $\alpha>0$ such that $\Delta(r)\geq r^\varepsilon$ for all $r\in ]0,\alpha[$.
Hence we obtain
\begin{eqnarray*}
	\Pi_1(t)&=&\frac{1}{(2\pi t)^{\frac{n}{2}}}\int_0^\alpha e^{-\frac{r^2}{2t}}r^{n-1}\Delta(r)\;dr\\
	&\geq&\frac{1}{(2\pi t)^{\frac{n}{2}}}\int_0^\alpha e^{-\frac{r^2}{2t}}r^{n-1+\varepsilon}\;dr\\
	&=&\frac{1}{(2\pi t)^{\frac{n}{2}}}\int_0^{\alpha/\sqrt{2t}} e^{-s^2}(s\sqrt{2t})^{n-1+\varepsilon}\sqrt{2t}\;ds\\
	&=&\pi^{-\frac{n}{2}}t^{\frac{\varepsilon}{2}}\int_0^{\alpha/\sqrt{2t}} e^{-s^2}s^{n-1+\varepsilon}\;ds.
\end{eqnarray*}
We have that 
$$\int_0^\infty e^{-s^2}s^{n-1+\varepsilon}\;ds>0.$$ 
Hence
$\Pi_1(t)\geq t^{\varepsilon}$ and therefore $\Phi_A(t)\geq t^\varepsilon$ for all sufficiently small $t$ and we are done.
\hfill$\Box$

\vs{0.5cm}
For $n\in\IN$ we denote by $\ma{S}_n$ the collection of all globally subanalytic subsets of $\IR^n$.
For $X\in \ma{S}_n$ let $F_X$ be the uniquely determined series in $\ma{A}_{n-2}$ such that $\Phi_X(t)\sim F_X(t)$ as $t\to 0$. By $F_X(0)$ we denote the constant term of the series $F_X$. Note that $\Phi_X$ and hence $F_X$ are invariant under applying an orthogonal transformation to $X$. Note also that $F_X(0)\in [0,1]$ since $\Phi_X(t)\in [0,1]$ for all $t\in\IR_{\geq 0}$.
We set 
$$\ma{F}:\ma{S}_n\to \ma{A}_{n-2}, X\mapsto F_X.$$
By $\ma{U}_{n,0}$ we denote the collection of all germs at $0$ of open globally subanalytic subsets of $\IR^n$.
By Theorem 3.6 we obtain a well-defined map
$$\mathfrak{F}:\ma{U}_{n,0}\to \ma{A}_{n-2}, U\mapsto \mathfrak{F}_U,$$
where $\mathfrak{F}_U:=F_V$ for some representative $V\in\ma{S}_n$ of $U$. 
Note that the map $\mathfrak{F}$ is not injective since it is for example invariant under orthogonal transformations.
But the empty germ is the only germ in $\ma{U}_{n,0}$ which is mapped to the zero series.

\vs{0.5cm}
We give a geometric criterion for the non-vanishing of the constant term of the asymptotic expansion which involves the notion of tangent cone (see Kurdyka and Raby [18]).

\vs{0.5cm}
{\bf 3.7 Theorem}

\vs{0.1cm}
{\it The following are equivalent:
\begin{itemize}
	\item[(i)] The constant term of the asymptotic expansion of $\Phi_A(t)$ at $0$ does not vanish.
	\item[(ii)] The tangent cone of $A$ at $0$ has full dimension. 
\end{itemize}}
{\bf Proof:}

\vs{0.1cm}
(i) $\Rightarrow$ (ii):
Assume that the tangent cone $T$ of $A$ at $0$ has dimension less than $n$. We find constants $r,c,\sigma\in \IR_{>0}$ with $\sigma>1$ such that
$$A\cap B(0,r)\subset \Big\{x\in \IR^n\;\Big\vert\;\mathrm{dist}(x,T)\leq c|x|^\sigma\Big\}=:C.$$
It suffices to show that the constant term $F_C(0)$ of $F_C$ vanishes.
Assume not. Then $F_C(0)>0$. Let $k\in \IN$ with $kF_C(0)>1$. For $s\in \IR_{>0}$ let $C(s):=C\cap B(0,s)$. By Theorem 3.6 we have
$F_C=F_{C(s)}$ for every $s>0$.
Since $\dim(T)<n$ and since $\sigma>1$
we find some $s>0$ and $k$ orthogonal $n\times n$ matrices $O_1,\ldots,O_k$ such that
$O_i\big(C(s)\big)\cap O_j\big(C(s)\big)=\{0\}$ for all distinct $i,j\in\{1,\ldots,k\}$. 
Let $D:=\bigcup_{i=1}^kO_i\big(C(s)\big)$. Since $F_D=\sum_{i=1}^k F_{C(s)}=kF_C$ we conclude that
$F_D(0)=kF_C(0)>1$, contradiction. 

\vs{0.2cm}
(ii) $\Rightarrow$ (i):  
Since the tangent cone of $A$ at $0$ has dimension $n$ we find some half line $l$ with start point $0$ and constants $c,r\in \IR_{>0}$ such that
$$E:=\Big\{x\in \IR^n\;\Big\vert\;\mathrm{dist}(x,l)\leq c|x|\Big\}\cap B(0,r)\subset A.$$
There is $k\in \IN$ and there are orthogonal $n\times n$-matrices $O_1,\ldots,O_k$ such that $\bigcup_{i=1}^kO_i(E)=\IR^n\cap B(0,r)$. 
We conclude that
$$F_A(0)\geq F_E(0)= \frac{1}{k}\sum_{i=1}^kF_{O_i(E)}(0)\geq \frac{1}{k}F_{\bigcup_{i=1}^kO_i(E)}(0)=\frac{1}{k}F_{\IR^n\cap B(0,r)}(0)=\frac{1}{k}F_{\IR^n}(0)=\frac{1}{k}.$$
\hfill$\Box$

\vs{0.5cm}
The next result shows that the asymptotic expansion at $0$ is in general not convergent. This implies that the function $\Phi_A$ is in general not definable in the o-minimal structure $\IR_{\an,\exp}$.

\vs{0.5cm}
{\bf 3.8 Example}

\vs{0.1cm}
{\it Consider the globally subanalytic subset $A$ of $\IR^2$ which is given in polar coordinates by
	$$\Big\{(r,\varphi)\in \IR_{>0}\times\; ]-\pi,\pi[\;\Big\vert\; r<\frac{1}{2}, 0<\varphi<\frac{r}{1-r^2}\Big\}.$$
	Then the following holds:
	\begin{itemize}
		\item[(1)] The asymptotic expansion of $\Phi_A$ at $0$ is not convergent.
	    \item[(2)] The function $\Phi_A$ is not definable in the o-minimal structure $\IR_{\an,\exp}$.
    \end{itemize}}

\vs{0.1cm}
{\bf Proof:}

\vs{0.1cm}
(1): Let
$\delta:\;]0,1/2[\;\to \IR, r\mapsto r/(1-r^2)$. 
We obtain for $r\in\; ]0,1/2[$ that
$$\Delta(r)=\int_{B_r}\;d\varphi=\int_{0}^{\delta(r)}\;d\varphi=
\delta(r).$$
Hence 
\begin{eqnarray*}
	\Phi_A(t)&=&\frac{1}{2\pi t}\int_0^{1/2}e^{-\frac{r^2}{2t}}r\delta(r)\;dr\\
	&=&\frac{1}{2\pi t}\int_0^{1/2}e^{-\frac{r^2}{2t}}\big(\sum_{k=1}^\infty r^{2k}\big)\;dr\\
	&=&\pi^{-1}\sum_{k=1}^\infty c_k\Big(\int_0^{1/(2\sqrt{2}t)}e^{-s^2}s^{2k}\;ds\Big)t^{k-\frac{1}{2}}
\end{eqnarray*}
where $c_k=2^{k-1/2}$. By the proof of Claim 1 in the proof of Theorem 3.5 we see that
$\Phi_A(t)\sim \sum_{k=1}^\infty E_kt^{k-1}$
where
$$E_k:=\frac{c_k}{\pi}\int_0^\infty e^{-s^2}s^{2k+1}\;ds.$$
By the substitution $\sigma=s^2$ we obtain 
$$\int_0^\infty e^{-s^2}s^{2k}\;ds=\frac{1}{2}\int_0^\infty e^{-\sigma}\sigma^{k-1}\;d\sigma=\frac{1}{2}\Gamma(k)=\frac{1}{2}(k-1)!$$
where $\Gamma$ is the Gamma function.
Hence 
$F_A=\sum_{k=1}^\infty E_kt^{k-1}\in \ma{A}_0\setminus \ma{A}_0^\mathrm{conv}$.

\vs{0.2cm}
(2): This follows from (2) and [9, Corollary 5.5].
\hfill$\Box$

\vs{0.5cm}
The logarithm may appear in the asymptotic expansion if $n\geq 3$ as the following example shows.

\vs{0.5cm}
{\bf 3.9 Example}

\vs{0.1cm}
{\it Consider the globally subanalytic subset $A$ of $\IR^3$ which is given in polar coordinates by
$$B:=\Big\{(r,\varphi,\vartheta)\in \IR_{>0}\times \;]-\pi,\pi[\;\times \;]-\pi/2,\pi/2[\;\Big\vert\; r<1, r<\varphi<1, 0<\vartheta<\arcsin\big(r/\varphi\big)\Big\}.$$
Then the asymptotic expansion of $\Phi_A$ at $0$ contains logarithmic terms.}

\vs{0.1cm}
{\bf Proof:}

\vs{0.1cm
}We obtain for $r\in \;]0,1[$ by Fubini's theorem that
$$\Delta(r)=\int_{B_r}C_3(\vartheta)\;d(\varphi,\vartheta)=\int_{\varphi=r}^{1}\int_{\vartheta=0}^{\arcsin(r/\varphi)}\cos(\vartheta)\;d\vartheta\;d\varphi=
-r\log r.$$
Hence 
\begin{eqnarray*}
	\Phi_A(t)
	&=&-\frac{1}{(2\pi t)^{\frac{3}{2}}}\int_0^{1}e^{-\frac{r^2}{2t}}r^3\log r\;dr\\
	&=&ct^{\frac{1}{2}}\int_{0}^{1/\sqrt{2t}}e^{-s^2}s^3\log(s\sqrt{2t})\;ds\\
	&=&ct^{\frac{1}{2}}\Big[\int_0^{1/\sqrt{2t}}e^{-s^2}s^3\log(\sqrt{2}s)\;ds\Big)+\frac{1}{2}\log t\Big(\int_0^{1/\sqrt{2t}}e^{-s^2}s^3\;ds\Big)\Big]
\end{eqnarray*}
where $c=\sqrt{2}/\pi^{\frac{3}{2}}$.
By the proof of Claim 1 in the proof of Theorem 3.5 we see that
$$\Phi_A(t)\sim (E_0+E_1\log t)t^{\frac{1}{2}}$$
	where
	$$E_0:=c\int_0^\infty e^{-s^2}s^3\log(\sqrt{2s})\;ds$$
	and
	$$E_1:=\frac{c}{2}\int_0^\infty e^{-s^2}s^3\;ds.$$
Since $E_1>0$ we obtain $F_A\in \ma{A}_1\setminus\ma{A}_0$.
\hfill$\Box$

\vs{0.5cm}
{\bf 3.10 Remark}

\vs{0.1cm}
Although $\Phi_A(t)$ has the Ilyashenko type of asymptotic expansion it is in general not in the Ilyashenko quasianalytic class (see Ilyashenko and Yakovenko [11, Section 24], [12, 13] or Speissegger [22] for the definitions). Moreover, although in the case $n=2$ no logarithmic terms occur, the function $\Phi_A(t)$ is in general not definable in the o-minimal structure $\IR_\ma{Q}$ which relies on this quasianalytic class (see [16]). This follows from [16, Theorem B] and the fact that $\Phi_A(t)$ also has in general exponentially small parts which are not dedected by the asymptotic expansions above. For example, from its asymptotic expansion (see [9, (5.7)]) we see that the error function is not definable in $\IR_\ma{Q}$.

\vs{1cm}
\noi \footnotesize{\centerline{\bf References}
	\begin{itemize}
		\item[(1)] 
		R. Cluckers, G. Comte, D. Miller, J.-P. Rolin, T. Servi:
		Integration of Oscillatory and Subanalytic Functions.
		arXiv 1601.01850.
		\item[(2)] 
		R. Cluckers, D. Miller:
		Stability under integration of sums of products of real globally subanalytic functions and their logarithms.
		{\it Duke Math. J.} {\bf 156} (2011), no. 2, 311-348.
		\item[(3)] 
		R. Cluckers, D. Miller:
		Loci of integrability, zero loci, and stability under integration for constructible functions on Euclidean space with Lebesgue measure.
		{\it Int. Math. Res. Not.} {\bf 2012}, no. 14, 3182-3191.
		\item[(4)]
		R. Cluckers, D. Miller:
		Lebesgue classes and preparation of real constructible functions.
		{\it J. Funct. Anal.} {\bf 264} (2013), no. 7, 1599-1642.
		\item[(5)] 
		G. Comte, J.-M. Lion, J.-P. Rolin:
		Nature log-analytique du volume des sous-analytiques.
		{\it Illinois J. Math.} {\bf 44} (2000), no. 4, 884-888.
		\item[(6)] 
		L. van den Dries:
		A generalization of the Tarski-Seidenberg Theorem, and some nondefinability results.
		{\it Bull. Amer. Math. Soc. (N.S.)} {\bf 15} (1986), no. 2, 189-193.
		\item[(7)] 
		L. van den Dries: 
	    Tame Topology and O-Minimal
		Structures. {\it London Math. Soc. Lecture Notes Series} {\bf
		248}, Cambridge University Press, 1998.
		\item[(8)] 
		L. van den Dries: O-minimal structures and real analytic geometry. In: B. Mazur et al. (Eds.), Current Developments in Mathematics, 1998. Proceedings of the Conference, Cambridge, MA, USA, 1998.
		International Press, 1999, pp. 105-152.
		\item[(9)]
		L. van den Dries, A. Macintyre, D. Marker:
		Logarithmic-exponential power series.
		{\it J. London Math. Soc.} (2) {\bf 56} (1997), no. 3, 417-434.
		\item[(10)]
		L. van den Dries, C. Miller:
		Geometric categories and o-minimal structures.
		{\it Duke Math. J.} {\bf 84} (1996), no. 2, 497-540.
		\item[(11)] 
		Y. Ilyashenko, S. Yakovenko: Lectures on analytic differential equations. American
		Mathematical Society, 2008.
		\item[(12)] 
		T. Kaiser: The Riemann mapping theorem for semianalytic domains and o-minimality.
		{\it Proc. Lond. Math. Soc. (3)} {\bf 98} (2009), no. 2, 427-444.
		\item[(13)] 
		T. Kaiser: The Dirichlet problem in the plane with
		semianalytic raw data, quasianalyticity and o-minimal structures.
		{\it Duke Math. Journal} {\bf 147}, no. 2 (2009), 285-314.
		\item[(14)] 
		T. Kaiser:
		First order tameness of measures.
		{\it Ann. Pure Appl. Logic} {\bf 163} (2012), no. 12, 1903-1927.
		\item[(15)] 
		T. Kaiser:
		Integration of semialgebraic functions and integrated Nash functions.
		{\it Math. Z.} {\bf 275} (2013), no. 1-2, 349-366.
		\item[(16)]
		T. Kaiser, J.-P. Rolin, P. Speissegger: Transition
		maps at non-resonant hyperbolic singularities are o-minimal. 
		{\it J. Reine Angew. Math.}
		{\bf 636} (2009), 1-45.
		\item[(17)] 
		I. Karatzas, S. E. Shreve: Brownian Motion and Stochastic Calculus. Spriner, 1996.
		\item[(18)] 
		K. Kurdyka, G. Raby:
		Densit\'e des ensembles sous-analytiques.
		{\it Ann. Inst. Fourier (Grenoble)} {\bf 39} (1989), no. 3, 753-771.
		\item[(19)]
		R. S. Lehman: Development of the mapping function at an analytic corner.
		{\it Pacific J. Math.} {\bf 7} (1957), no. 3, 1437-1449.
		\item[(20)] 
		J.-M. Lion and J.-P. Rolin:
		Int\'{e}gration des fonctions sous-analytiques et volumes des sous-ensembles sous-analytiques.
		{\it Ann. Inst. Fourier (Grenoble)} {\bf 48} (1998), no. 3, 755-767.
		\item[(21)]
		P. Speissegger: The Pfaffian closure of an o-minimal structure. 
		{\it J. Reine Angew. Math.} {\bf 508} (1999), 189-211.
		\item[(22)] 
		P. Speissegger: Quasianalytic Ilyashenko algebras.
		{\it Canadian Journal of Math.}, doi.org/10.4153/GJM-2016-048-x.
		\item[(23)] 
		W. Wasow: Asymptotic development of the solution of
		Dirichlet's problem at analytic corners. 
		{\it Duke Math. Journal} {\bf 24} (1957), 47-56.
\end{itemize}}

\vs{0.5cm}
Tobias Kaiser\\
University of Passau\\
Faculty of Computer Science and Mathematics\\
tobias.kaiser@uni-passau.de\\
D-94030 Germany

\vs{0.2cm}
Julia Ruppert\\
University of Passau\\
Faculty of Computer Science and Mathematics\\
julia.ruppert@uni-passau.de\\
D-94030 Germany}
\end{document}